\documentclass[3p]{elsarticle}

\usepackage{lineno,hyperref}

\usepackage{amsmath}
\usepackage{amsfonts}
\usepackage{amssymb}

\newtheorem{thm}{Theorem}

\newtheorem{cor}[thm]{Corollary}
\newproof{pf}{Proof}
\newdefinition{rmk}{Remark}

\newcommand{\dd}{\mathrm{d}}
\newcommand{\R}{\mathbb{R}}\newcommand{\C}{\mathbb{C}}

\usepackage{numcompress}

\begin{document}

\begin{frontmatter}
\title{Nondegeneracy of heteroclinic orbits for a class of potentials on the plane}


\author{Jacek Jendrej\corref{cor2}}
\ead{jendrej@math.univ-paris13.fr}
\address{CNRS \& LAGA, Universit\'e Sorbonne Paris Nord, 99 av Jean-Baptiste Cl\'ement, 93430 Villetaneuse, France}

\author{Panayotis Smyrnelis\corref{cor1}}
\ead{psmyrnelis@bcamath.org}
\cortext[cor1]{Corresponding author}
\address{Basque Center for Applied Mathematics, Alameda de Mazarredo 14, 48009 Bilbao, Spain}




\begin{abstract}
In the scalar case, the nondegeneracy of heteroclinic orbits is a well-known property, commonly used in problems involving nonlinear elliptic, parabolic or hyperbolic P.D.E. On the other hand, Schatzman proved that in the vector case this assumption is generic, in the sense that for any potential $W:\R^m\to\R$, $m\geq 2$, there exists an arbitrary small perturbation of $W$, such that for the new potential minimal heteroclinic orbits are nondegenerate. However, to the best of our knowledge, nontrivial explicit examples of such potentials are not available. In this paper, we prove the nondegeneracy of heteroclinic orbits for potentials $W:\R^2\to [0,\infty)$ that can be written as $W(z)=|f(z)|^2$, with $f:\C \to \C$ a holomorphic function.
\end{abstract}


\end{frontmatter}


\section{Introduction and main results}\label{sec:sec2}

Given a smooth nonnegative potential $W:\R^m\to [0,\infty)$ ($m\geq 1$) vanishing on a set $A$ of isolated points, i.e.
\begin{equation}\label{ww}
W\geq 0,\text{ and } W(u)=0 \Leftrightarrow u\in A,
\end{equation}
a \emph{heteroclinic orbit} (also called \emph{kink}) is a solution $e \in C^2(\R;\R^m)$ of the Hamiltonian system
\begin{equation}\label{system}
  e''(x)=\nabla W(e(x)), \ x\in \R,
\end{equation}
such that
\begin{equation}\label{limhet}
  \lim_{x\to -\infty}e(x)=a^- \text{ and } \lim_{x\to +\infty}e(x)=a^+ \text{ for some } a^-\neq a^+, \, \{a^-,a^+\}\subset A.
\end{equation}
That is, the heteroclinic orbit $e$ \emph{connects} at $\pm\infty$, two distinct global minima of $W$. We also recall that system \eqref{system} is the Euler-Lagrange equation of the energy functional
\begin{equation}\label{ene}
E(u):=\int_\R\Big[\frac{1}{2}|u'(x)|^2+W(u(x))\Big]\dd x, \  u \in H^1_{\mathrm{loc}}(\R;\R^m),
\end{equation}
and that the heteroclinic orbit $e$ is called \emph{minimal} if it is a minimizer of $E$ in the class $$\{u \in H^1_{\mathrm{loc}}(\R;\R^m):\lim_{x\to \pm\infty}u(x)=a^\pm\}.$$

The first existence proofs of a heteroclinic connection in the vector case were given by Rabinowitz \cite{rabinowitz} and by Sternberg \cite{sternberg1,sternberg2}. For more recent developments on the heteroclinic connection problem we refer to \cite{bon,indiana, antonop, monteil, sternberg3, fusco1,fusco2,book}. The aforementioned works provide sufficient conditions for the existence of heteroclinic orbits, in various settings. We also point out that if the minima $a^\pm$ of $W$ are assumed to be nondegenerate, in the sense that
\begin{equation}\label{nondeg}
\text{the Hessian matrices $D^2W(a^\pm)$ are positive definite},
\end{equation}
then \emph{every} heteroclinic orbit $e$ connecting $a^\pm$, approaches its limits at exponential rates, i.e.
\begin{equation}\label{expest}
\lim_{x\to\pm\infty} e^{k|x|}(|e(x)-a^\pm|+|e'(x)|)=0,
\end{equation}
holds for a constant $k>0$ independent of $e$ (cf. for instance \cite[Proposition 2.4]{book}).

Now, given a potential $W:\R^m\to[0,\infty)$, and a heteroclinic orbit $e$ connecting the nondegenerate zeros $a^\pm$, we consider the linearization operator of \eqref{system} at $e$, defined by
\begin{equation}
L:D(L)=H^2(\R;\R^m)\to L^2(\R;\R^m), \ Lh:=-h''+D^2W(e)h.
\end{equation}
It is easy to see that the operator $L$ is self-adjoint. Moreover, according to a theorem of Volpert et al. \cite{volp},
property \eqref{nondeg} implies that the essential spectrum of $L$ is $[M,+\infty)$, where $M>0$ is the minimum of the
lowest eigenvalues of $D^2W(a^-)$ and $D^2W(a^+)$. It is clear, by differentiating \eqref{system}, that $e'$ is an eigenvector of $L$ relative to the eigenvalue $0$. We shall say that a heteroclinic orbit $e$ is \emph{nondegenerate} if
\begin{equation}\label{kern}
\text{$L$ is nonnegative and } \dim \ker L=1.
\end{equation}

Note that if the heteroclinic orbit $e$ is minimal, then the first condition above
follows automatically, since the second variation of the energy satisfies
\begin{equation}
0 \leq \frac{\dd^2}{\dd \lambda^2}\Big|_{\lambda=0}E(e+\lambda h)=\int_\R(|h'|^2+h^\top D^2W(e)h)=\int_\R(-h''\cdot h+h^\top D^2W(e)h)=\langle Lh,h\rangle.
\end{equation}
Thus, for minimal $e$ the study of the nondegeneracy boils down to determining the kernel of $L$.

In the scalar case $m=1$, it is well-known that the heteroclinic orbit $e$ of a double well potential $W:\R\to[0,\infty)$ is unique (up to translations), minimal, and nondegenerate. The latter property is commonly used in problems involving nonlinear elliptic, parabolic or hyperbolic P.D.E. (cf. for instance \cite[Lemma 6.1]{annals} or \cite{henry,jendrej,kow1} for some relevant applications). On the other hand, to the best of our knowledge, nontrivial examples of potentials satisfying \eqref{kern} are not available in the vector case $m \geq 2$. However, Schatzman \cite[Section 4]{schatzman} proved that the nondegeneracy of \emph{minimal} heteroclinics is a \emph{generic} assumption, in the sense that
\begin{itemize}
\item[1)] if \eqref{kern} is satisfed for a given potential, it also holds for a small perturbation of that potential;
    \item[2)] given a potential $W$, and a minimal heteroclinic $e$ of \eqref{system}, there exists an arbitrarily small perturbation of $W$, such that for the new potential $e$ is still a minimal heteroclinic, and it satisfies \eqref{kern}.
\end{itemize}

The study of the nondegeneracy condition \eqref{kern} in \cite{schatzman} was motivated by the construction of heteroclinic double layers for the elliptic system $\Delta u(x)=\nabla W(u(x))$, $u:\R^2\to\R^2$. We point out that the existence of such a solution was initially established by Alama, Bronsard and Gui \cite{abg}, for symmetric potentials. Assumption \eqref{kern} is also relevant in the context of nonlinear evolution systems. We mention for instance the recent paper \cite{fusco3} by G. Fusco, where hypothesis \eqref{kern} is crucial to deduce the dynamics of multi-kinks under the parabolic flow $u_t(x,t)=u_{xx}(x,t)-\nabla W(u(x,t))$, $x\in (0,1)$, $t>0$.

In view of the aforementioned works, and the possible extension to the vector case of results established in the scalar case, it is important to determine potentials for which \eqref{kern} holds. In the present note, we provide a class of potentials defined on the plane, and satisfying \eqref{kern}. More precisely, we have the following result.

\begin{thm}\label{th1}
Let $W:\R^2\to [0,\infty)$ be a potential such that
\begin{equation}\label{hol}
W(z)=|f(z)|^2, \text{ where $f:\C\to\C$ is a holomorphic function. }
\end{equation}
Then, every heteroclinic orbit $e$ connecting two nondegenerate zeros $a^\pm$ of $W$, is nondegenerate (cf. \eqref{kern}).
\end{thm}
\begin{rmk}
In Theorem~\ref{th1}, we do not assume that $e$ is minimal.
\end{rmk}

We observe that potentials of the form \eqref{hol} include the natural class of potentials given as products of squares of distances from a finite set of points: $W(u)=\Pi_{i=1}^N |u-a_i|^2$. The proof of Theorem \ref{th1} is based on the complex analysis methods introduced by Alikakos, Betel\'u, and  Chen (cf. \cite{chen}), to solve the heteroclinic connection problem. For potentials $W$ satisfying the assumptions of Theorem \ref{th1}, it is established in \cite[Theorem 2]{chen} that there exists at most one heteroclinic orbit between any pair of wells of $W$. In addition (cf. \cite[Theorem 1]{chen}), if $g$ is a primitive of $f$, the image by $g$ of the heteroclinic orbit $x\mapsto e(x)$ is a line segment with end points $g(a^-)$ and $g(a^+)$.

As a consequence of Theorem \ref{th1}, we deduce, following \cite[Lemma 4.5]{schatzman}, a coercivity formula for the energy $E$, and for the second variation of $E$.

\begin{cor}\label{cor2}(cf. \cite[Lemma 4.5]{schatzman})
Under the assumptions of Theorem \ref{th1}, we set $e^T(x):=e(x-T)$, $\forall T\in\R$, and $\mathcal{C}(e):=\{e^T:T\in\R\}\subset e+H^1(\R;\R^2)$.
We also denote by $d$ the distance induced in the affine space $e+H^1(\R;\R^2)$ by the $H^1$ norm. Then, there exist constants $\alpha,\beta,\gamma,\delta>0$, such that
\begin{equation}\label{form2}
 \int_\R(|h'|^2+h^\top D^2W(e)h)\geq \alpha\|h\|_{H^1(\R;\R^2)}^2-\beta\frac{|\langle h,e'\rangle_{L^2(R;\R^2)}|^2}{\|e'\|_{L^2(\R;\R^2)}^2}, \  \forall h \in H^1(\R;\R^2),
\end{equation}
and
\begin{equation}\label{form1}
d(u,\mathcal{C}(e))\leq \gamma\Rightarrow E(u)-E(e)\geq \delta (d(u,\mathcal{C}(e)))^2.
\end{equation}
In particular, the heteroclinic orbit $e$ is a local minimum of $E$ in the affine space $e+H^1(\R;\R^2)$.
\end{cor}

\section{Proofs of Theorem \ref{th1} and Corollary \ref{cor2}}\label{sec:sec3}

\begin{pf}[Theorem \ref{th1}]

Proceeding as in \cite{chen}, we identify $e:=(e_1,e_2)\in \R^2$ with $e_1+i e_2 \in \C$, and notice that \eqref{system} is equivalent to
\begin{equation}\label{syshol}
e''=2f(e)\overline{f'(e)}.
\end{equation}
In addition, we have
\begin{equation}\label{equi}
\frac{1}{2}|e'|^2=f(e) \overline{f(e)},
\end{equation}
since the heteroclinic orbit $e$ satisfies the equipartition relation $\frac{1}{2}|e'|^2=W(e)$. Let $l$ be the arclength parameter defined by
$$l(x)=\int_{a^-}^x \sqrt{W( e(t))} \: |e'(t)| \dd t.$$
By computing
\begin{subequations}
\begin{equation}\label{subeq1}\frac{\dd}{\dd l}=(\sqrt{W(e)} \: |e'|)^{-1}\frac{\dd}{\dd x}=(\sqrt{W(e)}\sqrt{2W(e)} \: )^{-1}\frac{\dd}{\dd x}=\frac{1}{\sqrt{2}}(\overline{ f(e)} f(e))^{-1}  \frac{\dd}{\dd x},
\end{equation}
and
\begin{equation}\label{subeq2}
\frac{\dd g(e)}{\dd l}=\frac{1}{\sqrt{ 2}}\frac{g'(e)e'}{f(e) \overline{f(e)}}=\frac{1}{\sqrt{ 2}}\frac{f(e)e'}{f(e) \overline{f(e)}}=
\frac{1}{\sqrt{ 2}}\frac{e'}{\overline{f(e)}},
\end{equation}
we deduce that
\begin{equation}\label{subeq3}
\frac{\dd^2 g(e)}{\dd l^2}=\frac{\overline{f(e)} \: e''-e' \overline{f'(e) e'}}{2\overline{f(e)}\:^2 f(e) \overline{f(e)}}=\frac{e''-2f(e) \overline{f'(e)}}{ 2 \: f(e) \overline{f(e)}\:^2}=0,
\end{equation}
\end{subequations}
where both for \eqref{subeq1} and \eqref{subeq3} we used the equipartition relation \eqref{equi}. This proves that $\frac{\dd}{\dd l} g(e)=m$ is a constant. Moreover, it follows from \eqref{subeq2} and \eqref{equi} that
$$|m|^2=\Big|\frac{\dd}{\dd l} \: g(e)\Big|^2=\frac{1}{2} \frac{|e'|^2}{|f(e)|^2}=1,$$
thus $|m|=1$, and
\begin{equation}\label{dere}
e'=\sqrt{2} m \overline{f(e)}.
\end{equation}
Our next claim is that for every $h=(h_1,h_2) \in H^1(\R;\R^2)\sim H^1(\R;\C)$, we have
\begin{equation}\label{claim}
\int_\R(|h'|^2+h^\top D^2W(e)h)=\int_\R |h'-\sqrt{2}m\overline{f'(e)h}|^2.
\end{equation}
By expanding the right hand side of \eqref{claim}, we get
\begin{equation}\label{claim1}
\int_\R |h'-\sqrt{2}m\overline{f'(e)h}|^2=\int_\R(|h'|^2+2|f'(e)|^2|h|^2-\sqrt{2}\Re(\bar m f'(e)2hh')),
\end{equation}
and since $2hh'=(h^2)'$, an integration by parts gives
\begin{align}\label{claim2}
-\sqrt{2}\Re\Big(\int_\R(\bar m f'(e)2hh')\Big)&=\sqrt{2}\Re\Big(\bar m\int_\R( f''(e)e'h^2)\Big)\nonumber\\
&=2\Re\Big(\int_\R( f''(e)\overline{f(e)}h^2)\Big), \text{ in view of \eqref{dere}.}
\end{align}
In order to finish the proof of claim \eqref{claim}, we need to check that
\begin{equation}
h^\top D^2W(e)h=2|f'(e)|^2|h|^2+2\Re( f''(e)\overline{f(e)}h^2).
\end{equation}
This is an algebraic identity. We will check that $h^\top D^2W(z)h=2|f'(z)|^2|h|^2+2\Re( f''(z)\overline{f(z)}h^2)$ holds for all $h,z\in \C$. To this end, we expand $W(z+h)$ in powers of $h$, and examine the quadratic part:
$$W(z+h)=|f(z+h)|^2=|f(z)|^2+2\Re(f(z)\overline{f'(z)h})+|f'(z)|^2|h|^2+\Re( f''(z)\overline{f(z)}h^2)+o(|h|^2).$$
Now that \eqref{claim} is established, we are ready to determine the kernel of the operator $L$. For any $h \in \ker L$, it is clear  that
$$
\begin{gathered}
0=\langle Lh,h\rangle =  \int_\R(-h''\cdot h+h^\top D^2W(e)h)=\int_\R(|h'|^2+h^\top D^2W(e)h)
=\int_\R |h'-\sqrt{2}m\overline{f'(e)h}|^2 \\
\Rightarrow h'=\sqrt{2}m\overline{f'(e)h}.
\end{gathered}$$
On the other hand, equation $h'=\sqrt{2}m\overline{f'(e)h}$ can be written as a linear homogeneous system:
\begin{equation}
\begin{cases}
h'_1(x)=a_1(x)h_1(x)+a_2(x)h_2(x)\\
h'_2(x)=a_2(x)h_1(x)-a_1(x)h_2(x).
\end{cases}
\end{equation}
Recall that in view of \eqref{dere} $e'$ satisfies this system, hence the Wronskian $h_1e'_2-h_2e'_1=\mathrm{Const}$. Finally, since we assume decay as $x\to\pm\infty$ (cf. \eqref{expest}), this constant equals $0$, in other words:
\begin{equation*}
h(x)=\lambda(x)e'(x), \ \lambda(x)\in\R.
\end{equation*}
\end{pf}
But then, we have
\begin{align*}
h'(x)&=\lambda'(x)e'(x)+\lambda(x)e''(x)=\lambda'(x)e'(x)+\lambda(x)\sqrt{2}m\overline{f'(e(x))e'(x)}\\
&=\lambda'(x)e'(x)+\sqrt{2}m\overline{f'(e(x))h(x)},
\end{align*}
and this implies that $\lambda'\equiv 0$. Therefore, we obtain
\begin{equation*}
h=\lambda e', \ \lambda\in\R,
\end{equation*}
which completes the proof.\qed

\begin{pf}[Corollary \ref{cor2}]
In view of \eqref{kern}, we denote by $\lambda > 0$ the lower bound of the spectrum of $L$ without $0$. It follows that for any $h \in H^1(\R;\R^2)$, we have
\begin{equation}\label{spec}
\int_\R(|h'|^2+h^\top D^2W(e)h)\geq \lambda\|h\|_{L^2(\R;\R^2)}^2-\lambda\frac{|\langle h,e'\rangle_{L^2(R;\R^2)}|^2}{\|e'\|_{L^2(\R;\R^2)}^2}.
\end{equation}
On the other hand, it is clear that
\begin{equation}\label{spec2}
\int_\R(|h'|^2+h^\top D^2W(e)h)\geq \frac{1}{2}\|h\|_{H^1(\R;\R^2)}^2-\mu\|h\|_{L^2(\R;\R^2)}^2,
\end{equation}
holds for a constant $\mu>0$. Thus, we conclude that
\begin{equation}\label{spec3}
(1+\frac{\lambda}{\mu})\int_\R(|h'|^2+h^\top D^2W(e)h)\geq \frac{\lambda}{2\mu}\|h\|_{H^1(\R;\R^2)}^2-\lambda\frac{|\langle h,e'\rangle_{L^2(R;\R^2)}|^2}{\|e'\|_{L^2(\R;\R^2)}^2},
\end{equation}
which proves \eqref{form2}. The second coercivity formula \eqref{form1} is established in \cite[Lemma 4.5]{schatzman}. We point out that although \cite[Lemma 4.5]{schatzman} is stated for a \emph{minimal} nondegenerate heteroclinic orbit, it also holds under the weaker assumptions \eqref{kern}. Indeed, the arguments in the proof of \cite[Lemma 4.5]{schatzman} are based on properties \eqref{kern}, and not on minimality.
\qed
\end{pf}
\section*{Acknowledgments}
J. Jendrej was supported by ANR-18-CE40-0028 project ESSED.

P. Smyrnelis was supported by REA - Research Executive Agency - Marie Sk{\l}odowska-Curie Program (Individual Fellowship 2018) under Grant No. 832332, by the Basque Government through the BERC 2018-2021 program, by the Spanish State Research Agency through BCAM Severo Ochoa excellence accreditation SEV-2017-0718, and through the project PID2020-114189RB-I00 funded by Agencia Estatal de Investigaci\'on (PID2020-114189RB-I00 / AEI / 10.13039/501100011033).

\end{document}